\documentclass[12pt,twoside,leqno]{article}
\usepackage{latexsym}
\def\im{\mathop{\rm Im}\nolimits}
\def\dom{\mathop{\rm Dom}\nolimits}

\def\implies{\; \Longrightarrow \;}

\begin{document}
\begin{center}
{\bf  COMBINATORIAL RESULTS FOR CERTAIN SEMIGROUPS OF ORDER-DECREASING PARTIAL
ISOMETRIES OF A FINITE CHAIN}\\[4mm]
\textbf{F. Al-Kharousi, R. Kehinde and A. Umar}\\

\end{center}

\newtheorem{theorem}{{\bf Theorem}}[section]
\newtheorem{prop}[theorem]{{\bf Proposition}}
\newtheorem{lemma}[theorem]{{\bf Lemma}}
\newtheorem{corollary}[theorem]{{\bf Corollary}}
\newtheorem{remark}[theorem]{{\bf Remark}}
\newtheorem{conj}[theorem]{{\bf Conjecture}}
\newcommand{\pf}{\smallskip\noindent {\em Proof.}\ \  }
\newcommand{\qed}{\hfill $\Box$\medskip}

\newcommand{\inv}{^{-1}}

\begin{abstract}
Let ${\cal I}_n$ be the symmetric inverse semigroup on $X_n = \{1,
2, \ldots , n\}$ and let ${\cal DDP}_n$ and ${\cal ODDP}_n$ be its
subsemigroups of order-decreasing partial isometries and of
order-preserving and order-decreasing partial isometries of $X_n$,
respectively. In this paper we investigate the cardinalities of some
equivalences on ${\cal DDP}_n$ and ${\cal ODDP}_n$ which lead naturally
to obtaining the order of the semigroups. 
\footnote{Financial support from Sultan Qaboos University Internal Grant:
IG/SCI/DOMS/13/06 is gratefully acknowledged.}\end{abstract}


\section{Introduction and Preliminaries}

Let $X_n=\{1,2, \ldots, n\}$ and ${\cal I}_n$ be the partial
one-to-one transformation semigroup on $X_n$ under composition of
mappings. Then ${\cal I}_n$ is an {\em inverse} semigroup (that is,
for all $\alpha \in {\cal I}_n$ there exists a unique $\alpha' \in
{\cal I}_n$ such that $\alpha = \alpha\alpha'\alpha$ and $\alpha' =
\alpha'\alpha\alpha'$). The importance of ${\cal I}_n$ (more
commonly known as the {\em symmetric inverse semigroup or monoid}) to
inverse semigroup theory may be likened to that of the symmetric
group ${\cal S}_n$ to group theory. Every finite inverse semigroup
$S$ is embeddable in ${\cal I}_n$, the analogue of Cayley's theorem
for finite groups. Thus, just as the study of symmetric, alternating and
dihedral groups has made a significant contribution to group theory,
so has the study of various subsemigroups of ${\cal I}_n$, see for
example \cite{Bor, Fer1, Fer2, Gar, Lar, Uma1, Uma2}.

A transformation $\alpha \in {\cal I}_n$ is said to be
{\em order-preserving (order-reversing)} if $(\forall x,y \in \
\dom \alpha)\ x \leq y \implies x\alpha \leq y\alpha \ (x\alpha
\geq y\alpha)$ and, an {\em isometry (or distance-preserving)} if
($\forall x,y \in \dom \alpha) \mid x-y\mid = \mid x\alpha -y\alpha\mid$.
We shall denote by ${\cal DP}_n$ and ${\cal ODP}_n$, the semigroups
of partial isometries and of order-preserving partial isometries of
an $n-$chain, respectively. Eventhough semigroups of partial isometries on
more restrictive but richer mathematical structures have been studied
by Wallen \cite{Wal}, and Bracci and Picasso \cite{Bra} the study of the
corresponding semigroups on chains was only initiated recently by
Al-Kharousi et al. \cite{Kha1, Kha2}. A little while later, Kehinde et al.
\cite{Keh} studied ${\cal DDP}_n$ and ${\cal ODDP}_n,$ the order-decreasing
analogues of ${\cal DP}_n$ and ${\cal ODP}_n$, respectively.

Analogous to Al-Kharousi et al. \cite{Kha2}, this paper investigates the
combinatorial properties of ${\cal DDP}_n$ and ${\cal ODDP}_n$, thereby
complementing the results in Kehinde et al. \cite{Keh} which dealt
mainly with the algebraic and rank properties of these semigroups.
In this section we introduce basic definitions and terminology as well as
quote some elementary results from Section 1 of Al-Kharousi et al. \cite{Kha1}
and Kehinde et al. \cite{Keh} that will be needed in this paper. In Section 2
we obtain the cardinalities of two equivalences defined on ${\cal ODDP}_n$ and
${\cal DDP}_n$. These equivalences lead to formulae for the orders of ${\cal ODDP}_n$
and ${\cal DDP}_n$ as well as new triangles of numbers that were as a result of
this work recently recorded in \cite{Slo}.

For standard concepts in semigroup and symmetric inverse semigroup
theory, see for example \cite{How, Lip}. In particular $E(S)$
denotes the set of idempotents of $S$. Let
\begin{eqnarray} \label{eqn1.1} {\cal DDP}_n= \{\alpha \in {\cal DP}_n:
(\forall \ x\in\dom \alpha) \ x\alpha \leq x\}.
\end{eqnarray}
\noindent be the subsemigroup of ${\cal I}_n$ consisting of all
order-decreasing partial isometries of $X_n$. Also let
\begin{eqnarray} \label{eqn1.2} {\cal ODDP}_n= \{\alpha \in {\cal DDP}_n:
(\forall \ x,y\in \dom \alpha)\ x\leq y \Longrightarrow x\alpha\leq
y\alpha\}
\end{eqnarray}
\noindent be the subsemigroup of ${\cal DDP}_n$ consisting of all
order-preserving and order-decreasing partial isometries of $X_n$. Then
we have the following result.

\begin{lemma} \label{lem1.1} ${\cal DDP}_n$ and ${\cal ODDP}_n$ are
subsemigroups of ${\cal I}_n$.
\end{lemma}

\begin{remark}\label{rem1} ${\cal DDP}_n={\cal DP}_n \cap{\cal I}_n^-$ and
${\cal ODDP}_n={\cal ODP}_n \cap{\cal I}_n^-$, where ${\cal I}_n^-$
is the semigroup of partial one-to-one order-decreasing
transformations of $X_n$ \cite{Uma1}.
\end{remark}

Next, let $\alpha$ be an arbitrary element in ${\cal I}_n$. The {\em height} or
{\em rank} of $\alpha$ is $h(\alpha)= \mid \im \alpha\mid$, the {\em right [left]
waist} of $\alpha$ is $w^+(\alpha) = max(\im \alpha)\,
[w^-(\alpha) = min(\im \alpha)]$, the {\em right [left]
shoulder} of $\alpha$ is $\varpi^+(\alpha) = max(\dom \alpha)$\,
[$\varpi(\alpha) = min(\dom \alpha)]$, and {\em fix} of $\alpha$
is denoted by $f(\alpha)$, and defined by $f(\alpha)=|F(\alpha)|$, where
$$F(\alpha) = \{x \in X_n: x\alpha = x\}.$$

Next we quote some parts of \cite[Lemma 1.2]{Kha1} that will be needed
as well as state some additional observations that will help us understand more
the cycle structure of order-decreasing partial isometries.

\begin{lemma} \label{lem1.2} Let $\alpha\in {\cal DP}_n$. Then we have
the following:
\begin{itemize}
\item [(a)] The map $\alpha$ is either order-preserving or order-reversing. Equivalently,
$\alpha$ is either a translation or a reflection.
\item [(b)] If $f(\alpha)=p>1$ then
$f(\alpha)=h(\alpha)$. Equivalently, if $f(\alpha)>1$ then $\alpha$
is a partial identity.
\item [(c)] If $\alpha$ is order-preserving and $f(\alpha)\geq 1$ then $\alpha$
is a partial identity.
\item [(d)] If $\alpha$ is order-preserving then it is either strictly order-decreasing\\
($x\alpha < x$ for all $x$ in $\dom \alpha$) or  strictly order-increasing ($x\alpha > x$
for all $x$ in $\dom \alpha$) or a partial identity.
\item [(e)] If $F(\alpha)=\{i\}$ (for $1\leq i\leq n$) then for all $x\in \dom \alpha$ we have that
$x+x\alpha=2i$.
\item [(f)] If $\alpha$ is order-decreasing and $i\in F(\alpha)$ ($1\leq i\leq n)$ then for all $x\in \dom
\alpha$ such that $x < i$ we have $x\alpha =x$.
\item [(g)] If $\alpha$ is order-decreasing and $F(\alpha)=\{i\}$ then
$\dom \alpha\subseteq \{i, i+1,\ldots,n\}$.
\end{itemize}
\end{lemma}

\section{Combinatorial results}

For a nice survey article concerning combinatorial problems in the
symmetric inverse semigroup and some of its subsemigroups we refer
the reader to Umar \cite{Uma2}.

\noindent As in Umar \cite{Uma2}, for natural numbers $n\geq
p\geq m\geq 0$ and $n\geq i\geq 0$ we define

\begin{eqnarray} \label{eqn3.1} F(n;p_i)= \mid\{\alpha \in S:
h(\alpha)=\mid \im \alpha\mid = i \}\mid,
\end{eqnarray}

\begin{eqnarray} \label{eqn3.2} F(n;m_i)= \mid\{\alpha \in S:
f(\alpha)= i \}\mid
\end{eqnarray}

\noindent where $S$ is any subsemigroup of ${\cal I}_n$. From
\cite[Proposition 2.4]{Kha2} we have

\begin{theorem}\label{thrm2.1} Let $S = {\cal ODP}_n.$
Then $F(n;p)= \frac{(2n-p+1)}{p+1}{n\choose p}$, where $n\geq p \geq
1$.
\end{theorem}

We now have

\begin{prop}\label{prop2.2} Let $S = {\cal ODDP}_n$.
Then $F(n;p)= \pmatrix{n+1\cr p+1}$, where $n\geq p \geq 1$.
\end{prop}

\pf By virtue of Lemma \ref{lem1.2}[d] and Theorem \ref{thrm2.1} we see that

\begin{eqnarray*}
F(n;p) & = & \frac{1}{2}\left[\frac{2n-p+1}{p+1}{n\choose p}-{n\choose p}\right]+{n\choose p}\\
 & = &  \frac{1}{2}\left[\frac{2(n-p)}{p+1}{n\choose p}\right]+{n\choose p}\\
 & = & \frac{n-p}{p+1}{n\choose p}+{n\choose p}
  = {n\choose p+1} + {n\choose p} = {n+1\choose p+1}.\\
\end{eqnarray*}
\qed

The proof of the next lemma is routine using Proposition \ref{prop2.2}

\begin{lemma} \label{lem2.3} Let $S={\cal ODDP}_n$. Then
$F(n;p)=F(n-1;p-1)+F(n-1;p)$, for all $n\geq p\geq 2$.
\end{lemma}

\begin{theorem}\label{thrm2.4} $\mid {\cal ODDP}_n\mid = 2^{n+1}-(n+1).$
\end{theorem}

\pf It is enough to observe that $\mid {\cal ODDP}_n\mid=
\sum_{p=0}^{n}F(n;p)$.

\begin{lemma} \label{lem2.5} Let $S={\cal ODDP}_n$. Then
$F(n;m)={n\choose m}$, for all $n\geq m\geq 1$.
\end{lemma}

\pf It follows directly from Lemma \ref{lem1.2}[b,c] and
the fact that all idempotents are necessarily order-decreasing. \qed

\begin{prop}\label{prop2.6} Let \,$U_n=\{\alpha\in {\cal ODDP}_n: f(\alpha)=0\}$.
Then $\mid {U_n}\mid=\\ \mid {\cal ODDP}_{n-1}\mid$.
\end{prop}

\pf The proof is similar to that of \cite[Theorem 4.3]{Uma1}. \qed

\begin{remark} The triangles of numbers $F(n;p)$ and $F(n;m)$,  have as a result
of this work appeared in Sloane \cite{Slo} as [A184049] and [A184050], respectively.
\end{remark}

Now we turn our attention to counting order-reversing partial isometries. First recall
from \cite[Section3.2(c)]{Keh} that order-decreasing and order-reversing partial isometries
exist only for heights less than or equal to $n/2$. We now have

\begin{lemma}\label{lem2.8} Let $S={\cal DDP}^*_n$ be the set of order-reversing partial
isometries of $X_n$. Then $F(n;p_0)= 1$ and
$F(n;p_1)=\pmatrix{n+1\cr 2}$, for all $n\geq 1$.
\end{lemma}

\pf These follow from the simple observation that $$\{\alpha\in
{\cal ODDP}_n: h(\alpha)=0\,\,\mbox{or}\,\, 1\}=\{\alpha\in {\cal
DDP}^*_n: h(\alpha)=0\,\,\mbox{or}\,\, 1\}$$ \noindent and Proposition
\ref{prop2.2}.\qed

\begin{lemma} \label{lem2.10} Let $\alpha\in {\cal
DDP}^*_n$. Then for all $p\geq 1$ we have  \\ $F(2p+1,p+1)=1$ and $F(2p,p)=3$.
\end{lemma}

\pf {\bf (i)} By Lemma \ref{lem1.2}[f,g] we see that for $i\in \{0,1, \ldots, p\}$,
${p+1+i\choose p+1-i}$ is the unique order-reversing isometry of height $p+1;$
 and {\bf (ii)} for $i\in \{0,1, \ldots, p-1\}$, ${p+i\choose p-i}$,
${p+1+i\choose p-i}$ and ${p+1+i\choose p+1-i}$
are the only order-reversing isometries of height $p$.\qed

The following technical lemma will be useful later.

\begin{lemma} \label{lem2.11} Let $\alpha\in {\cal DDP}^*_n$. Suppose
$\varpi^+(\alpha)-r \in \dom \alpha$ and $\varpi^+(\alpha)-s\notin \dom \alpha$
for all $1\leq s <r$. Then $\varpi(\alpha)>r$.
\end{lemma}

\pf By order-reversing we see that $(\varpi^+(\alpha))\alpha=w^-(\alpha)$ and
$(\varpi(\alpha))\alpha=w^+(\alpha)$. Thus
$\varpi^+(\alpha)-r\geq \varpi(\alpha)\implies \varpi^+(\alpha)-\varpi(\alpha)
\geq r.$ So by isometry we have
$w^+(\alpha)-w^-(\alpha)=\varpi^+(\alpha)-\varpi(\alpha)\geq r
\implies w^+(\alpha)\geq w^-(\alpha)+r \implies w^+(\alpha) > r \implies
\varpi(\alpha) >r,$ as required.
\qed

\begin{lemma} \label{lem2.12} Let $S={\cal
DDP}^*_n$. Then
$F(n;p)=F(n-2;p-1)+F(n-2;p)$, for all $n\geq p\geq 2$.
\end{lemma}

\pf Let $\alpha \in {\cal
DDP}^*_n$ and $h(\alpha)=p$. Define $A=\{\alpha \in {\cal
DDP}^*_{n-2}: h(\alpha)=p\}$ and $B=\{\alpha \in {\cal DDP}^*_{n-2}:
h(\alpha)=p-1\}$. Clearly, $A \cap B=\emptyset$. Define a map $\theta :
 \{\alpha \in {\cal DDP}^*_n: h(\alpha)=p\} \rightarrow  A \cup B$ by
$(\alpha)\theta=\alpha'$ where

\noindent {\bf (i)} $x\alpha' =x\alpha\, (x\in \dom \alpha),$ if $\alpha\in A$.
It is clear that $\alpha'$ is an order-decreasing isometry and $h(\alpha)=p$;

\noindent {\bf (ii)} if $\{n-1, n\}\subseteq \dom \alpha\}$ and $\alpha\in B$, let
$\dom \alpha'=\{x-1: x\in \dom \alpha\,\,\mbox{and}\,\,x<n\}$
and $(x-1)\alpha'=x\alpha-1 \leq x-1$ and so $\alpha'$ is order-decreasing
and $h(\alpha)=p-1$;

\noindent {\bf (iii)} if $\{n-2, n-1\}\subseteq \dom \alpha\}$ and $\alpha\in B$, let
$\dom \alpha'=\{x-1: x\in \dom \alpha\,\,\mbox{and}\,\,x<n-1\}$
and $(x-1)\alpha'=x\alpha-1 \leq x-1$ and so $\alpha'$ is order-decreasing
and $h(\alpha)=p-1$;

\noindent {\bf (iv)} otherwise, if $\alpha\in B$, let
$\dom \alpha'=\{x-r: x\in \dom \alpha\,\,\mbox{and}\,\,x<\varpi^+(\alpha)\}$, where $r$
is such that $\varpi^+(\alpha)-r\in \dom \alpha$ and $\varpi^+(\alpha)-s\notin \dom \alpha$
for all $1\leq s < r$. Define $(x-r)\alpha'=x\alpha-r \leq x-r$ and so $\alpha'$ is
order-decreasing and Lemma \ref{lem2.11} ensures that $h(\alpha)=p-1$.

Moreover, in (ii) and (iii), we have $\mid(x-1)\alpha'-(y-1)\alpha'\mid
=\mid(x\alpha -1)-(y\alpha -1)\mid = \mid x\alpha -y\alpha\mid = \mid x-y\mid = \mid (x-1)-(y-1)\mid,$
and in (iv), we have $\mid(x-r)\alpha'-(y-r)\alpha'\mid
=\mid(x\alpha -r)-(y\alpha -r)\mid = \mid x\alpha -y\alpha\mid = \mid x-y\mid =\\ \mid (x-r)-(y-r)\mid.$
Hence $\alpha'$ is an isometry.

Also observe that in (ii), we have $\varpi^+(\alpha')=n-2$; in (iii) we have $\varpi^+(\alpha')=n-3$;
and in (iv) we have $\varpi^+(\alpha')<n-3$. These observations coupled with the definitions of $\alpha'$
ensures that $\theta$ is a bijection.

\noindent To show that $\theta$ is onto it is enough to note that we can in a symmetric manner
define $\theta^{-1}$ from  $ A \cup B  \rightarrow \{\alpha \in {\cal DDP}^*_n: h(\alpha)=p\}$.
This establishes the statement of the lemma. \qed

The next lemma which can be proved by induction, is necessary.

\begin{lemma}\label{lem2.13} Let $S = {\cal DDP}^*_n$. Then we have the following:
    $$\sum_{i\geq 0}{{n-1-2i\choose 2}}= \left\{
\begin{array}{ll}
\frac{(n+1)(n-1)(2n-3)}{24},\,&\,\mbox{if}\,\, $n$\,\,\mbox{is odd}; \\
\frac{n(n-2)(2n+1)}{24},\,&\,\mbox{if}\,\, $n$\,\,\mbox{is even}.
\end{array}
\right.$$
\end{lemma}

\begin{lemma}\label{lem2.14} Let $S = {\cal DDP}^*_n$. Then we have the following:
    $$F(n;p_2)= \left\{
\begin{array}{ll}
\frac{(n+1)(n-1)(2n-3)}{24},\,&\,\mbox{if}\,\, $n$\,\,\mbox{is odd}; \\
\frac{n(n-2)(2n+1)}{24},\,&\,\mbox{if}\,\, $n$\,\,\mbox{is even}.
\end{array}
\right.$$
 \end{lemma}

\pf By applying Lemmas \ref{lem2.8} and \ref{lem2.12} sucessively we get
\begin{eqnarray*} F(n;p_2)& = & F(n-2;p_1)+F(n-2;p_2)=  F(n-2;p_2)+{n-1\choose 2}\\
& = & F(n-4;p_2)+{n-3\choose 2}+{n-1\choose 2}\\
& = & F(n-6;p_2)+{n-5\choose 2}+{n-3\choose 2}+{n-1\choose 2}.\\
\end{eqnarray*}
By iteration the result follows from Lemma \ref{lem2.13}  and the facts that $F(2;p_2)=0$ and
$F(3;p_2)=1={2\choose 2}$.
\qed

\begin{prop}\label{prop2.15} Let $S = {\cal DDP}^*_n$. Then for all
$\lfloor (n+1)/2\rfloor\geq p\geq 1$, we have
$F(n;p)= \left\{
\begin{array}{ll}
\frac{(n+1)(n-1)(n-3)\cdots(n-2p+3)(2n-3p+3)}{2^p(p+1)!},\,&\,
\mbox{if}\,\, $n$\,\,\mbox{is odd}; \\
\frac{n(n-2)(n-4)\cdots(n-2p+2)(2n-p+3)}{2^p(p+1)!},\,&\,\mbox{if}\,\,
$n$\,\,\mbox{is even}.
\end{array}
\right.$.
\end{prop}

\pf (By Induction).

Basis Step: $F(n;p_1)=\pmatrix{n+1\cr 1+1}=\pmatrix{n+1\cr 2}$  is true by
Lemma \ref{lem2.8} and the observation made in its proof, while the formula
for $F(n;p_2)$ is true by Lemma \ref{lem2.14}.

Inductive Step: Suppose $F(m;p)$ is true for all $\lfloor (n+1)/2\rfloor > m\geq p\geq 1$.

{\bf Case 1.} If $m$ is odd, consider (using the induction hypothesis)
\begin{eqnarray*} F(m+2;p)& = & F(m;p)+F(m;p-1)\\
& = & \frac{(m+1)(m-1)(m-3)\cdots(m-2p+3)(2m-3p+3)}{2^p(p+1)!}\\
& + & \frac{(m+1)(m-1)(m-3)\cdots(m-2p+5)(2m-3p+6)}{2^{p-1}p!}\\
& = & \frac{(m+3)(m+1)(m-1)\cdots(m-2p+5)(2m-3p+7)}{2^p(p+1)!},\\
\end{eqnarray*}
\noindent which is the formula for $F(m+2;p)$ when $m$ is odd.

{\bf Case 2.} If $m$ is even, consider (using the induction hypothesis)
\begin{eqnarray*} F(m+2;p)& = & F(m;p)+F(m;p-1)\\
& = & \frac{m(m-2)(m-4)\cdots(m-2p+2)(2m-p+3)}{2^p(p+1)!}\\
& + & \frac{m(m-2)(m-4)\cdots(m-2p+4)(2m-p+4)}{2^{p-1}p!}\\
& = & \frac{(m+2)m(m-2)\cdots(m-2p+4)(2m-p+7)}{2^p(p+1)!},\\
\end{eqnarray*}
\noindent which is the formula for $F(m+2;p)$ when $m$ is even.
\qed

\begin{prop}\label{prop2.16} Let $S = {\cal DDP}^*_n$ and let $b_n = \sum_{p\geq 0}{F(n;p)}$.
Then for $n\geq 0$, we have
\begin{enumerate}
\item $b_{2n+1}=5\cdot2^{n+1}-4n-8$;
\item $b_{2n}=7\cdot2^n-4n-6$.
\end{enumerate}
\end{prop}

\pf Apply induction and use the fact that $\mid {\cal DDP}^*_n\mid=
\sum_{p=0}^{n}F(n;p)$.

\begin{prop}\label{prop2.17} Let $S = {\cal DDP}_n$. Then

\begin{itemize}
\item[(1)] if $n$ is odd and $(n+1)/2\geq p\geq 2$ \\
$F(n;p)=\frac{(n+1)(n-1)(n-3)\cdots(n-2p+3)(2n-3p+3)}{2^p(p+1)!}+{n+1\choose p+1}$;
\item[(2)] if $n$ is even and $n/2\geq p\geq 2$ \\
$F(n;p)=\frac{n(n-2)(n-4)\cdots(n-2p+2)(2n-p+3)}{2^p(p+1)!}+{n+1\choose p+1}$;
\item[(3)] if $\lfloor(n+1)/2\rfloor < p$, $F(n;p)={n+1\choose p+1}$.
\end{itemize}
\end{prop}

\pf It follows from Propositions \ref{prop2.2} \& \ref{prop2.15} and
Lemmas \ref{lem1.2}[c] \& \ref{lem2.8}.\qed

Combining Theorem \ref{thrm2.4}, Lemmas \ref{lem1.2}[a,c] \& \ref{lem2.10}, Proposition
\ref{prop2.16} and the observation made in the proof of Lemma \ref{lem2.8} we get the
order of ${\cal DDP}_n$ which we record as a theorem below.

\begin{theorem}\label{thrm2.18} Let ${\cal DDP}_n$.
Then for all $n\geq 0$ we have
\begin{itemize}
\item[(1)] $\mid {\cal DDP}_{2n+1}\mid = 2^{2n+2}+5\cdot2^{n+1}-(2n^2+9n+12)$;
\item[(2)] $\mid {\cal DDP}_{2n}\mid = 2^{2n+1}+7\cdot2^n-(2n^2+7n+8)$.
\end{itemize}
\end{theorem}

\begin{lemma}\label{lem2.19} Let $S={\cal DDP}_n$. Then
$F(n;m)={n\choose m}$, for all $n\geq m\geq 2$.
\end{lemma}

\pf It follows directly from \cite[Lemma 3.18]{Keh} and the fact
that all idempotents are necessarily order-decreasing. \qed

\begin{prop}\label{prop2.20} Let $S={\cal DDP}_n$. Then $F(2n;m_1)= 2^{n+1}-2$
and $F(2n-1;m_1)= 3\cdot2^{n-1}-2$, for all $n\geq 1$.
\end{prop}

\pf Let $F(\alpha)=\{i\}$. Then by Lemma \ref{lem1.2}[e], for any $x\in
\dom \alpha$ we have $x+x\alpha=2i$. Thus, by Lemma
\ref{lem1.2}[g], there $2i-2$ possible elements for $\dom \alpha:
(x,x\alpha)\in \{(i,i),(i+1,i-1), (i+2,i-2), \ldots, (2i-1,1)\}$.
However, (excluding $(i,i)$) we see that there are
$\sum_{j=0}{i-1\choose j}=2^{i-1}$, possible partial isometries with
$F(\alpha)=\{i\}$, where $2i-1\leq n \Longleftrightarrow i\leq
(n+1)/2$. Moreover, by symmetry we see that $F(\alpha)=\{i\}$ and
$F(\alpha)=\{n-i+1\}$ give rise to equal number of decreasing
partial isometries. Note that if $n$ is odd (even) the equation $i=n-i+1$
has one (no) solution. Hence, if $n=2a-1$ we have
$$2\sum_{i=1}^{a-1}2^{i-1}+2^{a-1}=2(2^{a-1}-1)+ 2^{a-1}=3.2^{a-1}-2$$
\noindent decreasing partial isometries with exactly one fixed
point; if $n=2a$ we have

$$2\sum_{i=1}^{a}2^{i-1}=2(2^a-1)=2^{a+1}-2$$
\noindent decreasing partial isometries with exactly one fixed
point.\qed

\begin{theorem}\label{thrm2.21} Let ${\cal DDP}_n$.
Then $$a_n =\mid {\cal DDP}_n\mid =
3a_{n-1}-2a_{n-2}-2^{\lfloor\frac{n}{2}\rfloor}+{n+1},$$
\noindent with $a_0 = 1$ and $a_{-1} = 0$.
\end{theorem}

\pf It follows from Propositions \ref{prop2.6} \& \ref{prop2.20}, Lemma \ref{lem2.19}
and the fact that $\mid {\cal DDP}_n\mid= \sum_{m=0}^{n}F(n;m)$.\qed

\begin{remark} The triangle of numbers $F(n;m)$ and sequence $\mid{\cal DDP}_n\mid$
have as a result of this work appeared in Sloane \cite{Slo} as [A184051] and [A184052],
respectively. However, the triangles of numbers $F(n;p)$ for ${\cal DDP}_n$ and
${\cal DDP^*}_n$ and the sequence $\mid{\cal DDP}^*_n\mid$ are as at the time
of submitting this paper not in Sloane \cite{Slo}.  For some computed values of
$F(n;p)$, see Tables 3.1 and 3.2.
\end{remark}

\begin{center}

$$\begin{array}{|c|c|c|c|c|c|c|c|c|c|}
\hline
 \,\,\,\,\,n{\backslash}p&0&1&2&3&4&5&6&7&\sum F(n;p)=\mid {\cal DDP}^*_n \mid
\\ \hline 0&1&&&&&&&&1
 \\ \hline 1&1&1&&&&&&&2
 \\ \hline 2&1&3&0&&&&&&4
\\
\hline 3&1&6&1&0&&&&&8
\\
\hline 4&1&10&3&0&0&&&&14
\\
\hline 5&1&15&7&1&0&0&&&24
\\
\hline 6&1&21&13&3&0&0&0&&38
\\
\hline 7&1&28&22&8&1&0&0&0&60
\\
\hline
\end{array}$$
\end{center}

\begin{center}
Table 3.1
\end{center}

\begin{center}

$$\begin{array}{|c|c|c|c|c|c|c|c|c|c|}
\hline
 \,\,\,\,\,n{\backslash}p&0&1&2&3&4&5&6&7&\sum F(n;p)=\mid {\cal DDP}_n \mid
\\ \hline 0&1&&&&&&&&1
 \\ \hline 1&1&1&&&&&&&2
 \\ \hline 2&1&3&1&&&&&&5
\\
\hline 3&1&6&5&1&&&&&13
\\
\hline 4&1&10&13&5&1&&&&30
\\
\hline 5&1&15&27&16&6&1&&&66
\\
\hline 6&1&21&48&38&21&7&1&&137
\\
\hline 7&1&28&78&78&57&28&8&1&279
\\
\hline
\end{array}$$
\end{center}

\begin{center}
Table 3.2
\end{center}

\section{Number of ${\cal D}^*$-classes}

For the definitions of the Green's relations (${\cal L}, {\cal R}$ and ${\cal D}$) and their
starred analogues (${\cal L}^*, {\cal R}^*$ and ${\cal D}^*$), we refer the reader to
Howie \cite{How} and Fountain \cite{Fou2}, (respectively) or Ganyushkin and Mazorchuk \cite{Gan}.

First, notice that from \cite[Lemma 2.1]{Kha1} we deduce that number of ${\cal L}$-classes in
$K(n,p)=\{\alpha\in {\cal DP}_n: h(\alpha)=p\}$ (as well as the number of ${\cal R}$-classes there)
is ${n\choose p}$. To describe the ${\cal D}$-classes in  ${\cal DP}_n$ and  ${\cal ODP}_n$,
first we recall (from \cite{Kha1}) that the {\em gap} and {\em reverse gap of the image set of
$\alpha$ (with $h(\alpha)=p$) are ordered $(p-1)$-tuples} defined as follows:
$$g(\im \alpha)=(\mid a_2\alpha -a_1\alpha\mid, \mid a_3\alpha -a_2\alpha\mid, \ldots,
\mid a_p\alpha -a_{p-1}\alpha\mid)$$
\noindent and
$$g^R(\im \alpha)=(\mid a_p\alpha -a_{p-1}\alpha\mid),\ldots, \mid a_3\alpha -a_2\alpha\mid,
\mid a_2\alpha -a_1\alpha\mid),$$

\noindent where $\alpha=\pmatrix{a_1&a_2&\cdots&a_p\cr
a_1\alpha&a_2\alpha&\cdots&a_p\alpha}$ with $1\leq a_1<a_2<\cdots <a_p\leq n.$
Further, let $d_i=\mid a_{i+1}\alpha -a_i\alpha\mid$ for $i=1,2, \ldots, p-1$. Then
$$g(\im \alpha)=(d_1, d_2, \ldots,d_{p-1})\,\, \mbox{and} \,\,g^R(\im \alpha)=(d_{p-1}, d_{p-2},\ldots, d_1).$$
For example, if
$$\alpha=\pmatrix{1&2&4&7&8\cr
3&4&6&9&10},\beta=\pmatrix{2&4&7&8\cr
10&8&5&4}\in {\cal DP}_{10}$$

\noindent then $g(\im \alpha)=(1,2,3,1),\, g(\im \beta)=(2,3,1),\, g^R(\im \alpha)=(1,3,2,1)$
and $g^R(\im \beta)=(1,3,2).$ Next, let $d(n,p)$ be the number of distinct ordered $p$-tuples:
$(d_1, d_2, \ldots,d_p)$ with $\sum_{i=1}^{p}d_i=n$. This is clearly the number of
{\em compositions} of $n$ into $p$ parts. Thus, we have

\begin{lemma} \cite[p.151]{Rio} \label{lem3.1} $d(n,p)={n-1\choose p-1}.$
\end{lemma}

We shall henceforth use the following well-known binomial identity when needed:

$$\sum_{m=p}^{n}{m\choose p} = {n+1\choose p+1}.$$

We take this opportunity to state and prove a result which was omitted in \cite{Kha2}.

 \begin{theorem}\label{thrm3.1} Let $S = {\cal ODP}_n$. Then
 \begin{itemize}
 \item[(1)] the number of ${\cal D}$-classes in $K(n,p)\,(p\geq 1)$ is ${n-1\choose p-1}$;
 \item[(2)] the number of ${\cal D}$-classes in $S$ is $1+2^{n-1}$.
 \end{itemize}
   \end{theorem}

 \pf
 \begin{itemize}
 \item[(1)] It follows from \cite[Theorem 2.5]{Kha1}:
 $(\alpha, \beta)\in {\cal D}$ if and only if $g(\im \alpha)=g(\im \beta)$;
 \cite[Lemma 3.3]{Kha1}: $p-1\leq \sum_{i=1}^{p-1} d_i\leq n-1$; Lemma\,\ref{lem3.1};
 and so the number of ${\cal D}$-classes is
 $\sum_{i=p-1}^{n-1} d(i, p-1)= \sum_{i=p-1}^{n-1}{i-1\choose p-2}={n-1\choose p-1}.$
 \item[(2)] The number of ${\cal D}$-classes in $S$ is
 $1+\sum_{p=1}^{n}{n-1\choose p-1}=1+2^{n-1}$.
 \end{itemize}
 \qed

 The following results from \cite{Keh} will be needed:

 \begin{lemma} \label{lem3.2} \cite[Lemma 2.3]{Keh} Let $\alpha ,\beta \in {\cal DDP}_n$ or
${\cal ODDP}_n.$ \ Then
\begin{itemize}
\item[(1)] $\alpha \leq _{{\cal R^*}}\beta $ \,  if and only if\,
$\dom \alpha\subseteq \dom \beta$;
\item[(2)] $\alpha \leq _{{\cal L^*}}\beta $\,   if and only if  $\im \alpha
\subseteq \im \beta$;
\item[(3)] $\alpha \leq _{{\cal H^*}}\beta $ \,  if and only if\,
$\dom \alpha\subseteq \dom \beta$ and $\im \alpha \subseteq
\im \beta$.
\end{itemize}
\end{lemma}

\noindent From \cite[(3)]{Keh}, for $\alpha, \beta\in {\cal DDP}_n$, we have
$(\alpha, \beta)\in {\cal D^*}$ if and only if
\begin{eqnarray} \label{eqn3.2}
g(\im \alpha)= \left\{
\begin{array}{ll}
g(\im \beta);\,\, \mbox{or} \\
g^R(\im \beta),\,\,\mbox{if}\,\,p\leq a_p-a_1\leq (n-1)/2.
\end{array}
\right.
\end{eqnarray}

\noindent Similarly, from \cite[(4)]{Keh}, for $\alpha, \beta\in {\cal ODDP}_n$,
we have

\begin{eqnarray} \label{eqn3.2} (\alpha, \beta)\in {\cal D^*} \,\,\mbox{if and only if}\,\,
g(\im \alpha)=g(\im \beta).
\end{eqnarray}

Now a corollary of Theorem\,\ref{thrm3.1} follows:

\begin{corollary}\label{cor3.3} Let $S = {\cal ODDP}_n$. Then
 \begin{itemize}
 \item[(1)] the number of ${\cal D}^*$-classes in $K(n,p)\,(p\geq 1)$ is ${n-1\choose p-1}$;
 \item[(2)] the number of ${\cal D}^*$-classes in $S$ is $1+2^{n-1}$.
 \end{itemize}
   \end{corollary}

Observe that for all $\alpha \in {\cal DP}_n$ with $h(\alpha)=p$,

\begin{eqnarray}  \label{eqn3.3} a_p - a_1 = \sum_{i=1}^{p-1}(a_{i+1}-a_i)=\sum_{i=1}^{p-1}d_i,
\end{eqnarray}

\noindent where $g(\dom \alpha)= (d_1, d_2, \ldots, d_{p-1}).$  Moreover, an ordered $p$-tuple:
$(d_1, d_2, \ldots,d_p)$ is said to be {\em symmetric} if
$$(d_1, d_2, \ldots,d_p)=(d_1, d_2, \ldots,d_p)^R=(d_p, d_{p-1}, \ldots,d_1).$$
Now, let $d_s(n,p)$ be the number of distinct symmetric ordered $p$-tuples: \\
$(d_1, d_2, \ldots,d_p)$ with $\sum_{i=1}^{p}d_i=n$. Then we have

\begin{lemma} \cite[Lemma 3.5]{Kha2}\label{lem3.5} $d_s(n;p)= \left\{
\begin{array}{ll}
0,\,&\,\,\mbox{if}\,\, $n$\,\,\mbox{is odd and}\,\, $p$\,\,\mbox{is even}; \\
{\lfloor {\frac{n-1}{2}}\rfloor\choose \lfloor {\frac{p-1}{2}}\rfloor},\,&\,\,\mbox{otherwise}.
\end{array}
\right. $
\end{lemma}

Now by virtue of $(5)$ and \cite[Theorem 2.5]{Kha1}, it is not difficult to see that the
number of ${\cal D}^*$-classes in ${\cal DDP}_n$ is the same as the number of
${\cal D}$-classes in ${\cal ODP}_n$ less those pairs that are merged into single
${\cal D}^*$-classes in ${\cal DDP}_n$. Thus, we have

\begin{lemma} \label{lem3.6} Let $g(m,p)$ be the number of ${\cal D}$-classes in ${\cal ODP}_n$
(consisting of maps of height $p$ and $\sum d_i = m$) that are merged into single
${\cal D}^*$-classes in ${\cal DDP}_n$. Then $m\leq (n-1)/2$, and \\
$g(m,p)= \left\{
\begin{array}{ll}
\frac{1}{2}{m-1\choose p-2},\,&\,\,\mbox{if}\,\, $n$\,\,\mbox{is odd and}\,\, $p$\,\,\mbox{is odd}; \\
\frac{1}{2}[{m-1\choose p-2}-{\lfloor {\frac{m-1}{2}}\rfloor\choose \lfloor {\frac{p-2}{2}}\rfloor}],\,&\,\,\mbox{otherwise}.
\end{array}
\right. $
\end{lemma}

\pf The result follows from $(5)$, Lemmas \ref{lem3.1} \& \ref{lem3.5} and the
observation that\\
\begin{eqnarray*}
g(n,p) = \frac{d(n-1,p-1) - d_s(n-1,p-1)}{2}.
\end{eqnarray*}\qed

Now have the main result of this section.

\begin{theorem} \label{thrm3.7} Let $B(n,p)$ be the number of ${\cal D}$-classes in ${\cal ODP}_n$
(consisting of maps of height $p$) that are merged into single ${\cal D}^*$-classes in ${\cal DDP}_n$.
Then for $n\geq p\geq 1$, we have \\
$B(n,p)= \left\{
\begin{array}{ll}

\frac{1}{2}[{\lfloor{\frac{n-1}{2}}\rfloor\choose p-1}-
{\lfloor{\frac{n-1}{4}}\rfloor\choose \frac{p-1}{2}}],
\,&\,\,\mbox{if}\,\, $p$\,\,\mbox{is odd}; \\
\frac{1}{2}[{\lfloor{\frac{n-1}{2}}\rfloor\choose p-1}-
2{\lfloor{\frac{n-1}{4}}\rfloor\choose \frac{p}{2}}],
\,&\,\,\mbox{if}\,\,n\equiv 1,2\, (mod\,4),
\,\&\,$p$\,\,\mbox{is even}; \\
\frac{1}{2}[{\lfloor{\frac{n-1}{2}}\rfloor\choose p-1}-
2{\lfloor{\frac{n-3}{4}}\rfloor\choose \frac{p}{2}}-
{\lfloor{\frac{n-3}{4}}\rfloor\choose \frac{p-2}{2}}],
\,&\,\,\mbox{if}\,\,n\equiv -1,0\, (mod\,4),\,\&\,$p$\,\,\mbox{is even}.
\end{array}
\right. $
\end{theorem}

\pf The result follows from $(5)$, $(7)$ and Lemma \ref{lem3.6}. To see this, let
$n\equiv 0\,(mod\,4)$ and $p$ be even. Then $n=4k$ for some integer $k$, and \\
\begin{eqnarray*} B(n,p)& = & \sum_{m=p}^{\lfloor{\frac{n-1}{2}}\rfloor}g(m,p)
= \sum_{m=p}^{2k-1}g(m,p)\\
& = & g(p,p)+g(p+2,p)+\cdots +g(2k-2,p) \\
& + & g(p+1,p)+g(p+3,p)+\cdots +g(2k-1,p)\\
& = & \frac{1}{2}\left[{p-1\choose p-2}-{\frac{p-2}{2}\choose \frac{p-2}{2}}
+{p+1\choose p-2}-{\frac{p}{2}\choose \frac{p-2}{2}}
+\cdots +{2k-3\choose p-2} -{k-2\choose \frac{p-2}{2}}\right]\\
& + & \frac{1}{2}\left[{p\choose p-2}-{\frac{p}{2}\choose \frac{p-2}{2}}
+{p+2\choose p-2}-{\frac{p+2}{2}\choose \frac{p-2}{2}}
+\cdots +{2k-2\choose p-2} -{k-1\choose \frac{p-2}{2}}\right]\\
& = & \frac{1}{2}\left[{2k-1\choose p-1}-2{k-1\choose \frac{p}{2}}
-{k-1\choose \frac{p-2}{2}}\right]\\
& = & \frac{1}{2}\left[{\frac{n-2}{2}\choose p-1}-2{\frac{n-4}{4}\choose \frac{p}{2}}
-{\frac{n-4}{4}\choose \frac{p-2}{2}}\right].
\end{eqnarray*}
\noindent All the other cases are handled similarly.\qed

Now have the main result of this section.

\begin{corollary} \label{cor3.8} The number of ${\cal D}^*$-classes in ${\cal DDP}_n$
(consisting of maps of height $p\geq 1$) is ${n-1\choose p-1}- B(n,p)$.
\end{corollary}

\pf The result follows from Theorem \ref{thrm3.7} and the remarks preceding Lemma\,\ref{lem3.6}.\qed

\begin{corollary} The number of ${\cal D}^*$-classes in ${\cal DDP}_n$ denoted by $d_n$
is\\
 $d_n= \left\{
\begin{array}{ll}
2^{n-1}-2^{\lfloor{\frac{n-3}{2}}\rfloor}+\cdot2^{\lfloor{\frac{n+1}{4}}\rfloor},
\,&\,\,\mbox{if}\,\,  n\equiv -1,0\,(mod\,4); \\
2^{n-1}-2^{\lfloor{\frac{n-3}{2}}\rfloor}+3\cdot2^{\lfloor{\frac{n-3}{4}}\rfloor},
\,&\,\,\mbox{if}\,\,  n\equiv 1,2\,(mod\,4).
\end{array}
\right. $
\end{corollary}

\pf The result follows from Theorem \ref{thrm3.7} and Corollary \ref{cor3.8}. To see this, let
$n\equiv 1,2\,(mod\,4)$. Then $n=4k+1,4k+2$ for some integer $k$, and \\
\begin{eqnarray*} d_n& = & 1+\sum_{p=1}^{n}{n-1\choose p-1}-
\sum_{p=1}^{\lfloor{\frac{n-1}{2}}\rfloor}B(n,p)= 1+2^{n-1}-\sum_{p=1}^{2k}B(n,p)\\
& = & 1+2^{n-1}-[B(n,1)+B(n,3)+\cdots +B(n,2k-1)] \\
& - & [B(n,2)+B(n,4)+\cdots +B(n,2k)]\\
& = & 1+2^{n-1}-\frac{1}{2}\left[{2k\choose 0}-{k\choose 0}
+{2k\choose 2}-{k\choose 1}
+\cdots +{2k\choose 2k-2} -{k\choose k-1}\right]\\
& - & \frac{1}{2}\left[{2k\choose 1}-2{k\choose 1}
+{2k\choose 3}-2{k\choose 2}
+\cdots +{2k\choose 2k-1} -2{k\choose k}\right]\\
& = &  1+2^{n-1}-\frac{1}{2}\left[(2^{2k}-1)-(3\cdot2^k+1)+2\right]\\
& = &  2^{n-1}-2^{\lfloor{\frac{n-3}{2}}\rfloor}+3\cdot2^{\lfloor{\frac{n-3}{4}}\rfloor}.
\end{eqnarray*}
\noindent The case $n\equiv -1,0\,(mod\,4)$ is handled similarly.\qed

\noindent {\bf Acknowledgements}. The second named author would like
to thank Bowen University, Iwo and Sultan Qaboos University for
their financial support and hospitality, respectively.

\vspace{1cm}

\begin{center}
F. Al-Kharousi\\
Department of Mathematics and Statistics\\
Sultan Qaboos University \\
Al-Khod, PC 123 -- OMAN\\
E-mail:{\tt fatma9@squ.edu.om}

\end{center}

\begin{center}
R.\ Kehinde\\
Department of Mathematics and Statistics\\
Bowen University \\
P. M. B. 284, Iwo, Osun State\\
Nigeria.\\
E-mail:{\tt kennyrot2000@yahoo.com}
\end{center}

\begin{center}
 A.\ Umar\\
Department of Mathematics\\
Petroleum Institute, P. O. Box 2533 \\
Abu Dhabi, U. A. E.\\
E-mail:{\tt aumar@pi.ac.ae}

\end{center}

\end{document}